\documentstyle[11pt]{article}
\parskip =\medskipamount
\def\be{\begin{equation}}
\def\ee{\end{equation}}
\def\disp{\displaystyle}

\newcounter{fig}

\def\th{\theta}
\def\lam{\lambda}
\def\R{\rm {\sf R\hspace{-3.05mm}I\hspace{2.1mm}}}
\def\Z{\rm {\sf Z\hspace{-3.15mm}Z\hspace{0.6mm}}}

\begin{document}

\begin{titlepage}

\begin{center}
{\LARGE On the limiting power of set of knots generated by $1+1$-- and 
$2+1$-- braids

}
\vspace{0.4in}

{\large R Bikbov$^{\dag}$ \hspace{0.5cm} and \hspace{0.5cm}
S Nechaev$^{\ddag\dag}$} \bigskip

$^{\dag}$ {\sl L D Landau Institute for Theoretical Physics, \\ 117940, 
Moscow, Russia}

$^{\ddag}$ {\sl Institut de Physique Nucl\'eaire, Division de Physique
Th\'eorique$^{*}$, \\ 91406 Orsay Cedex, France}

\end{center}

\vspace{0.5in}

\begin{abstract}
We estimate from above the set of knots, $\Omega(n,\mu)$, generated by closure of 
$n$--string $1+1$-- and $2+1$--dimensional braids of irreducible length $\mu$
($\mu\gg 1$) in the limit $n\gg 1$.  
\end{abstract} 
\vspace{0.4in}

\vspace{2in}


\vspace{0.3in}

\noindent {\bf PACS:} 
\vspace{0.5in}

\hrule \footnotesize
$^{*}$ Unit\'e de Recherche des Universit\'es Paris XI et Paris VI
associ\'ee au C.N.R.S.

\end{titlepage}

\section{Introduction}

Besides the traditional fundamental topological issues concerning the
construction of new topological invariants, investigation of homotopic
classes and fibre bundles we mark a set of ajoint but much less studied
problems. First of all, we mean the problem of so-called ``knot entropy"
calculation. Most generally it can be formulated as follows. Take the
lattice ${\Z}^3$ embedded in the space ${\R}^3$. Let $\Omega$ be the
ensemble of all possible closed nonselfintersecting $N$--step paths with 
one common fixed point on ${\Z}^3$; by $\omega_N$ we denote the particular
configuration of the trajectory. The main question is: what is the fraction 
$P_N$ of the the trajectories $\omega_N\in\Omega$ belonging to some specific 
homotopic class characterized by the topological invarant ${\mbox{\sl Inv}}$ 
(we do not specify the way of defining the topological invariant). The 
distribution function $P\{\mbox{\sl Inv}\}$ satisfies the obvious 
normalization condition $\disp \sum_{{\rm all}\;\omega_N\in\Omega} 
P_N\{\mbox{\sl Inv}\}=1$.

In the present paper we pay attention to the statistical problem concerning 
the estimation of the set $\Omega=\{\Omega^{(1)}, \Omega^{(2)}\}$ of 
knots generated by closure of braids embedded in $1+1$-- and $2+1$-- 
dimensions (see the definitions below). 

The paper is organized as follows. Below we give the basic definitions of the 
standard $1+1$--dimensional and $2+1$--dimensional braid groups as well as 
formulate the basic results; the Section 2 is devoted to the estimations of 
the sets $\Omega^{(1)}$ and $\Omega^{(2)}$ using the concept of $1+1$-- 
and $2+1$-- locally--free groups; while in Conclusion we discuss in more 
details the corollaries following from our consideration.

\subsection{The basic definitions}

\noindent 1. The $1+1$--dimensional (``standard") braid group 
$B^{(1)}_{n+1}$ of $n+1$ strings has $n$ generators 
$\{\sigma_1,\sigma_2,\ldots, \sigma_n\; \mbox{and their inverses}\}$ 
(see fig.\ref{fig1}a) with the following relations:  
\be \label{1} 
\left\{\begin{array}{ll}  
\sigma_i\sigma_{i+1}\sigma_i = 
\sigma_{i+1}\sigma_i\sigma_{i+1} & \qquad (1\le i<n) \medskip \\ 
\sigma_i\sigma_j=\sigma_j\sigma_i & \qquad (|i-j|\ge 2) \medskip \\ 
\sigma_i\sigma_i^{-1}=\sigma_i^{-1}\sigma_i=e & 
\end{array} \right.
\ee

\noindent 2. The $2+1$--dimensional (``surface") braid group 
$B^{(2)}_{n+1}$ can be defined in the following way (see, for instance 
\cite{pr,birm}). Consider the two--dimensional lattice ${\Z}^2$ and take 
distinct points $P_1, P_2, \ldots, P_{n+1}\in {\Z}^2$. A $2+1$--braid of 
$n+1$ strings on ${\Z}^2$ based at $\{P_1, P_2, \ldots, P_{n+1}\}$ is an 
$n+1$--tuple $b=(b_1,\ldots, b_{n+1})$ of paths, $b_i: [1,N]\to {\Z}^2$, 
such that 

(i) $b_i(1)=P_i$ and $b_1(1)\in \{P_1, P_2, \ldots, P_{n+1}\}$ $\forall 
i\in \{1,\ldots n+1\}$;

(ii) $b_i(t)\neq b_j(t)$ $\forall \{i,j\}\in \{1,\ldots n+1\}$, $i\neq j$; 
$t\in [1,N]$.

\noindent The braid group $B^{(2)}_{n+1}$ on ${\Z}^2$ based at $\{P_1, 
P_2, \ldots, P_{n+1}\}$ is the group of homotopy classes of braids based at 
$\{P_1, P_2, \ldots, P_{n+1}\}$. The group $B^{(2)}_{n+1}$ has $2(n\times 
n)$ generators $\{(\sigma_{11}^{(x)},\sigma_{11}^{(y)}),\ldots,\break 
(\sigma_{1n}^{(x)}, \sigma_{1n}^{(y)}); \ldots;(\sigma_{n1}^{(x)},
\sigma_{n1}^{(y)}),\ldots,(\sigma_{nn}^{(x)},\sigma_{nn}^{(y)})\; 
\mbox{and their inverses}\}$ (see fig.\ref{fig1}b) with the following 
relations:  
\be \label{1a} 
\left\{\begin{array}{ll}  
\sigma_{i,j}^{(x)}\sigma_{i+1,j}^{(x)}\sigma_{i,j}^{(x)} = 
\sigma_{i+1,j}^{(x)}\sigma_{i,j}^{(x)}\sigma_{i+1,j}^{(x)} & \qquad (1\le 
\{i,j\}\le n) \medskip \\
\sigma_{i,j}^{(y)}\sigma_{i,j+1}^{(y)}\sigma_{i,j}^{(y)} = 
\sigma_{i,j+1}^{(y)}\sigma_{i,j}^{(y)}\sigma_{i,j+1}^{(y)} & \qquad (1\le 
\{i,j\}\le n) \medskip \\
\sigma_{i,j}^{(x)}\sigma_{i,j}^{(y)}\sigma_{i,j}^{(x)} = 
\sigma_{i,j}^{(y)}\sigma_{i,j}^{(x)}\sigma_{i,j}^{(y)} & \qquad (1\le 
\{i,j\}\le n) \medskip \\
\sigma_{i_1,j_1}^{(x)}\sigma_{i_2,j_2}^{(x)}= 
\sigma_{i_2,j_2}^{(x)}\sigma_{i_1,j_1}^{(x)} & \qquad
(|i_1-i_2|>1\quad {\rm or}\quad |j_1-j_2|>0) 
\medskip \\ 
\sigma_{i_1,j_1}^{(x)}\sigma_{i_2,j_2}^{(y)}= 
\sigma_{i_2,j_2}^{(y)}\sigma_{i_1,j_1}^{(x)} & \qquad (i_2-i_1 \ne
\{0,1\}\quad {\rm or}\quad j_1-j_2 \ne \{0,1\}) \medskip \\ 
\sigma_{i,j}^{(x)}\left(\sigma_{i,j}^{(x)})\right)^{-1}= 
\sigma_{i,j}^{(y)}\left(\sigma_{i,j}^{(y)})\right)^{-1}=e &
\end{array} \right.
\ee

The braid groups $B^{(1)}_n$ and $B^{(2)}_n$ have the following general
properties:

-- Any arbitrary word written in terms of ``letters"---generators
of the groups $B^{(1)}_n$ or $B^{(2)}_n$---gives a particular braid. 

-- The {\it length}, $N$, of the braid is the total number of used
letters, while the {\it minimal irreducible length}, $\mu$, hereafter 
referred to as the ``primitive length" is the shortest noncontractible 
length of a particular braid which remains after applying of all possible 
group relations. Diagramatically the braid can be represented as a set of 
crossed strings going from the top to the bottom appeared after subsequent 
gluing the braid generators.

-- The closed braid appears after gluing the ``upper" and the ``lower"
free ends of the braid on the cylinder.

\subsection{The main results}

Our basic results might be formulated in a geometrically clear way.
Consider two sets of braids $\{B^{(1)}_n\}$ and $\{B^{(2)}_n\}$, 
embedded in $1+1$-- and $2+1$-- dimensions correspondingly. Let each 
particular braid has the primitive length $\mu$ and is represented by 
$n$ strings. 

Then:  
\begin{itemize} 
\item The set $\Omega^{(1)}(n,\mu)$ of knots which can be generated by 
the standard braids of given irreducible length $\mu$ ($\mu\gg 1$) from the set 
$\{B^{(1)}_n\}$ ($n={\rm const}\gg 1$) is restricted from above by the 
value 
\be \label{2}
\Omega^{(1)}(n,\mu)< \frac{32\pi^2}{\ln^42}\;\frac{2^n}{n^3}\;7^{\mu-1} 
\ee 
\item The set $\Omega^{(2)}(n,\mu)$ of knots which can be generated by 
the surface braids of given irreducible length $\mu$ ($\mu\gg 1$) from the set 
$\{B^{(2)}_n\}$ ($n={\rm const}\gg 1$) is restricted from above by the 
value 
\be \label{3}
\Omega^{(2)}(n,\mu)<\frac{32n^2}{\pi^2}\left(\frac{2n}{\ln n}\right)^{\mu-1}  
\ee
\end{itemize}
(See the Conclusion for more detailed discussion of the results (\ref{2}) and 
(\ref{3})).

\section{Combinatorics of words}

Any braid corresponds to some knot or link. The correspondence between 
braids and knots is not mutually single valued and each knot or link can be 
represented by infinite series of different braids. However, we can 
estimate from above the partition functions $\Omega^{(1)}(n,\mu)$ and 
$\Omega^{(2)}(n,\mu)$ of all possible knots generated by the ensemble of 
all $1+1$-- and $2+1$-- braids of primitive length $\mu$ using the 
following obvious fact. {\it The sets $\Omega^{(1)}(n,\mu)$ and 
$\Omega^{(2)}(n,\mu)$ are bounded from above by the number of all 
distinct words of the primitive length $\mu$ in $1+1$-- and $2+1$-- braid 
groups correspondingly.} Thus in what follows we are aimed in the 
estimation of the number of nonequivalent words in the standard and surface 
braid groups.

\subsection{Definitions of $1+1$-- (``standard") and $2+1$-- (``surface") 
locally free groups}

\noindent 1. Following the ideas of A.M. Vershik concerning the 
notion of the "local groups" \cite{ver1} and the papers \cite{negrve}, where 
the concept of a "locally free" group was proposed at first in the topological 
context, let us define the group, ${\cal LF}_{n+1}^{(1)}$, which has $n$ 
generators $\{f_1,\ldots,f_n\; \mbox{and their inverses}\}$ with the relations:
\be \label{4}
\left\{\begin{array}{ll} 
f_j f_k=f_k f_j \qquad \mbox{for} \quad |j-k|\ge 2 \medskip \\
f_i f_i^{-1}=e
\end{array}\right.
\ee
We call the group with relations (\ref{4}) the $1+1$--dimensional "locally 
free group", because each pair of generators $(f_j, f_{j\pm 1})$ produces a 
free subgroup of the group ${\cal LF}_{n+1}^{(1)}$.

The group ${\cal LF}_{n+1}^{(1)}$ can be obtained from the braid group 
$B_{n+1}^{(1)}$ if we replace the braiding ("Yang-Baxter-type") relations 
by the free ones. The geometrical interpretation of the generators of a 
group ${\cal LF}_{n+1}^{(1)}$ is shown in fig.\ref{fig2}a.

Apparently, in mathematical literature the notion similar to our "locally 
free group" appeared firstly in the paper \cite{cafo} devoted to the 
investigation of the combinatorial properties of rearrangements of 
sequences, known also as "partially commutative monoids" (see \cite{vien} 
and references therein).

\noindent 2. The $2+1$--dimensional (``surface") locally free group 
${\cal LF}^{(2)}_{n+1}$ has $2(n\times n)$ generators 
$\{(f_{11}^{(x)},f_{11}^{(y)}),\ldots,(f_{1n}^{(x)},f_{1n}^{(y)});
\ldots;(f_{n1}^{(x)},f_{n1}^{(y)}),\ldots,(f_{nn}^{(x)},f_{nn}^{(y)})\;
\mbox{and their inverses}\}$ with the following relations:  
\be \label{4a} 
\left\{\begin{array}{ll}  
f_{i_1,j_1}^{(x)}f_{i_2,j_2}^{(x)}= 
f_{i_2,j_2}^{(x)}f_{i_1,j_1}^{(x)} & \qquad (|j_1-j_2|>0\quad 
{\rm or}\quad |i_1-i_2|>1)
\medskip \\ 
f_{i_1,j_1}^{(x)}f_{i_2,j_2}^{(y)}= 
f_{i_2,j_2}^{(y)}f_{i_1,j_1}^{(x)} & \qquad (i_2-i_1 \ne \{0,1\}\quad 
{\rm or}\quad j_1-j_2\ne \{0,1\}) \medskip \\ 
f_{i,j}^{(x)}\left(f_{i,j}^{(x)})\right)^{-1}= 
f_{i,j}^{(y)}\left(f_{i,j}^{(y)})\right)^{-1}=e &
\end{array} \right.
\ee
Thus, we can construct the $2+1$--locally free group ${\cal 
LF}^{(2)}_{n+1}$ from the surface braid group $B^{(2)}_{n+1}$ if we 
replace the braiding relations of the neighbouring generators by the 
"full monodromy", i.e. by the free group relations---see the fig.\ref{fig2}b.

The following important properties of $1+1$-- and $2+1$-- locally free 
groups should be mentioned:
\begin{itemize}
\item[(i)] By definition the locally free groups ${\cal LF}^{(1)}_{n+1}$
and ${\cal LF}^{(2)}_{n+1}$ have less relations than the braid groups 
$B^{(1)}_{n+1}$ and $B^{(2)}_{n+1}$ correspondingly. 
{\it Thus, the number of distinct words of the primitive length $\mu$ in 
the $1+1$-- and $2+1$-- braid groups is bounded from above by the number 
of distinct words of the primitive length $\mu$ in the $1+1$-- and $2+1$-- 
locally free groups.}
\item[(ii)] By construction (compare figures \ref{fig1} and \ref{fig2}) the 
monodromy generators $f_i$ ($i\in[1,n]$) of the group ${\cal
LF}^{(1)}_{n+1}$ 
and $f^{(x,y)}_{i,j}$ ($\{i,j\}\in[1,n]$) of the group ${\cal 
LF}^{(2)}_{n+1}$ can be written as $f_i=\left(\sigma_i\right)^2$ 
($i\in[1,n]$) and $f^{(x,y)}_{i,j}=\left(\sigma^{(x,y)}_{i,j}\right)^2$ 
($\{i,j\}\in[1,n]$), where $\sigma_i$ and $\sigma^{(x,y)}_{i,j}$ are the 
generators of the groups $B^{(1)}_{n+1}$ and $B^{(2)}_{n+1}$ correspondingly. 
{\it Thus, the number of distinct words of the primitive length $2\mu$ in 
the $1+1$-- and $2+1$-- braid groups is bounded from below by the number of 
distinct words of the primitive length $\mu$ in the $1+1$-- and $2+1$-- 
locally free groups.}
\end{itemize}

\subsection{Computation of number of nonequivalent words in $1+1$-- and 
$2+1$-- locally free groups}

We derive explicitly the expressions of the numbers $V^{(1)}(n,\mu)$ and 
$V^{(2)}(n,\mu)$ of all nonequivalent primitive words of length $\mu$ in the 
groups ${\cal LF}_{n+1}^{(1)}$ and ${\cal LF}_{n+1}^{(2)}$ 
respectively. Our computations are based on the so-called "normal 
order" representation of words proposed by A.M. Vershik in \cite{ver2} 
(see also \cite{negrve}). 
\bigskip

\noindent {\sc The group ${\cal LF}_{n+1}^{(1)}$}. Let us represent each 
word $W_p$ of irreducible length $\mu$ in the group ${\cal LF}_{n+1}^{(1)}$ in 
the "standard" form 
\be \label{2:norm}
W_p=\left(f_{\alpha_1}\right)^{m_1}\left(f_{\alpha_2}\right)^{m_2}\ldots
\left(f_{\alpha_s}\right)^{m_s}
\ee
where $\sum_{i=1}^s |m_i|=\mu\; (m_i\neq 0\; \forall\; i;\; 1\le s\le
\mu$) and the sequence of generators $f_{\alpha_i}$ in Eq.(\ref{2:norm})
{\it for all distinct} $f_{{\alpha}_i}$ satisfies the following local rules
\cite{negrve} ("normal order" representation):
\begin{itemize}
\item[(i)] If $f_{\alpha_i}=f_1$, then $f_{\alpha_{i+1}}=f_2$;
\item[(ii)] If $f_{\alpha_i}=f_k$ ($2\le k\le n-1$), then
$f_{\alpha_{i+1}} 
\in \left\{f_1,\ldots,f_{k-1},f_{k+1},\right\}$; 
\item[(iii)] If $f_{\alpha_i}=f_n$, then $f_{\alpha_{i+1}}\in
\left\{f_1,\ldots,f_{n-1}\right\}$.  
\end{itemize}

The rules (i)--(iii) give the prescription how to encode and enumerate all
distinct primitive words in the group ${\cal LF}_{n+1}^{(1)}$. If the
sequence of generators in the primitive word $W_p$ does not satisfy the
rules (i)-(iii), we commute the generators in the word $W_p$ until the
normal order is restored.  Hence, the normal order representation enables
one to give the unique coding of all nonequivalent primitive words in our 
group.

Let $\th_n(m)$ be the number of all distinct sequenses of $m+1$
generators, $1\le m\le \mu-1$, satisfying the rules (i), (ii), (iii). The 
calculation of the number of distinct primitive words $V^{(1)}(n,\mu)$ of 
given primitive length $\mu$ is now straightforward:
\be \label{v1}
V^{(1)}(n,\mu)=\sum_{m=1}^{\mu-1}2^{m+1}\left(\mu-1 \atop m\right)\th_n(m).
\ee
The combinatorial factor $2^{m+1}\left(\mu-1 \atop m\right)$ in Eq.(\ref{v1}) 
is the number of all primitive words of length $\mu$ written in a normal order 
form for the {\it fixed sequence} of $m+1$ generators. 

Our approach to the computation of $\th_n(m)$ is based on the consideration 
of a  "correlation function" $\th_n(x,x_0,m)$ which is defined as the number 
of all distinct sequences of $m+1$ generators satisfying the rules (i), (ii), 
(iii), beginning with the generator $f_{x_0}$ and ending with the generator 
$f_x$. It is easy to write an evolution equation for $\th(x,m) \equiv 
\th_n(x,x_0,m)$ with the "time" $m$: 
\be \label{5}
\th(x,m+1)=\th(x-1,m)+\sum_{y=x+1}^{n}\th(y,m) 
\ee
This equation should be completed by initial and boundary conditions 
\be \label{6}
\left\{\begin{array}{l}
\th(x,0)=\delta_{x,x_0} \\ \th(0,m)=\th(n+1,m)=0 
\end{array}\right.
\ee

We solve the boundary problem (\ref{5})--(\ref{6}) in the limit $n\gg 1$ 
supposing the periodical boundary conditions on the segment $[0,n+1]$. 
Namely, we have:
\be \label{theta}
\left\{\begin{array}{l}
\disp \theta(x+1,m+1)-\theta(x,m+1)= 
\theta(x,m)-\theta(x-1,m)-\theta(x+1,m) \medskip \\
\disp \theta(x,0)=\delta_{x_0,x} \medskip \\
\disp \theta(0,m)=\theta(n+1,m)=0. 
\end{array}\right.
\ee
The substitution
\be 
\th(x,m)=\sum_{k=1}^{n} A_k \lambda_k^m \alpha_k(x)
\ee
enables us to pass to the following recursion relations:
\be \label{alpha} 
\left\{\begin{array}{l}
(\lambda_k+1)\alpha_k(x+1)-(\lambda_k+1)\alpha_k(x)+\alpha_k(x-1)=0 \\
\alpha_k(0)=\alpha_k(n+1)=0 \\
\end{array}\right. 
\ee
One can readily find the eigen--values and eigen--functions of (\ref{alpha}):
\be \label{eigen1d}
\begin{array}{l}
\disp \lam_k=4\cos^2\frac{\pi k}{n+1} - 1 \medskip \\ 
\disp \alpha_k(x)=
\frac{\sin\frac{\pi kx}{n+1}}{(2\cos\frac{\pi k}{n+1})^x},
\qquad k=1,\ldots,n \\ 
\end{array}
\ee 

As the function (\ref{alpha}) is not symmetric, the set of eigen--functions 
is not orthogonal on the segment $[0,n+1]$ and it is difficult to ensure the 
initial condition. It is convinient to pass from (\ref{alpha}) to symmetric 
problem. Consider a generating function
\be  
Z(x,s)=\sum_{m=0}^{\infty} s^m \th(x,m)
\ee
The equation for the function $Z(x,s)$ reads
\be 
\left\{\begin{array}{l}
(s+1)Z(x+1,s)-(s+1)Z(x,s) + sZ(x-1,s)=
\delta_{x,x_0-1}-\delta_{x,x_0}\\
Z(0,s)=Z(n+1,s)=0
\end{array}\right.
\ee
The last equation can be symmetrized via the substitution  
$Z(x,s)=A^x\varphi(x,s)$, where $A=\sqrt{\frac{s}{s+1}}$. Thus, we get
\be \label{7}
\left\{\begin{array}{l}
\varphi(x+1,s)-\frac{1}{A}\varphi(x,s)+\varphi(x-1,s)= \disp
\frac{A^{-x}}{\sqrt{s(s+1)}} (\delta_{x,x_0-1}-\delta_{x,x_0}) \medskip \\
\varphi(0,s)=\varphi(n+1,s)=0.
\end{array}\right.
\ee

Making use of the sin-Fourier transform, $\disp f(k,s)=\sum_{x=1}^{n} 
\varphi(x,s)\sin\frac{\pi kx}{n+1}$, let us rewrite (\ref{7}) in the
form
$$
\left(2\cos\frac{\pi k}{n+1} -\frac{1}{A}\right)f(k,s)=
\frac{A^{-x_0}}{\sqrt{s(s+1)}}\left(A\sin\frac{\pi k(x_0-1)}
{n+1}-\sin\frac{\pi kx_0}{n+1}\right). 
$$
The final explicit expression of the function $Z(x,s)$ reads as follows
$$ 
Z(x,s)=\frac{2}{(n+1)(s+1)}
\left(\frac{s}{s+1}\right)^\frac{x-x_0}{2}\sum_{k=1}^{n}
\frac{ \sqrt{\frac{s+1}{s}}\sin\frac{\pi kx_0}{n+1} -\sin\frac{\pi
k(x_0-1)}{n+1}}{\sqrt{ \frac{s+1}{s}}-2\cos\frac{\pi k}{n+1}}\sin
\frac{\pi kx}{n+1}.
$$
Now we can restore the function $\th(x,m)$ via contour integration
$$ 
\th(x,m)=\frac{1}{2\pi i}\oint\limits_C 
\frac{Z(x,s)}{s^{m+1}}ds, 
$$
where the contour $C$ surrounds the point $s=0$ and is displaced 
in the regularity area of the function $Z(s)\equiv Z(x,s)$. 
Hence 
$$
\th(x,m)=-\sum\limits_{s_k} \mbox{Res}
\left(\frac{Z(x,s_k)}{s_k^{m+1}}\right), 
$$ 
where $s_k$ are the poles out of the regularity area:  
$$
s_k=\frac{1}{4\cos^2\frac{\pi k}{n+1}-1} 
$$ 
(compare to (\ref{eigen1d})). We are interested only in the asymptotic 
behavior $m\gg 1$ of the function $\th_n(x,x_0,m)$ which is determined by 
the poles nearest to the origin, $s_1=s_{n-1}=\frac{1}{3}$ for $n \gg 1$. 
So we get 
\be 
\th_n(x,x_0,m)=\frac{4}{n+1}\sin\frac{\pi (x_0+1)}{n+1}
\sin\frac{\pi x}{n+1}\;2^{x_0-x}\;3^m 
\ee 
To find the function $\th_n(m)$ we should sum up $\th(x,x_0,m)$ over all 
$x$ and $x_0$: $\disp \th_n(m)=\break \sum_{x,x_0=1}^{n}\th_n(x,x_0,m)$. 
We obttain in the limit for $n={\rm const}\gg 1$ the following expression: 
\be \label{18}
\th_n(m)=\frac{16\pi^2}{\ln^4 2}\;\frac{2^n}{n^3}\;3^m
\ee
The whole number of nonequivalent words follows from Eq.(\ref{v1}) in 
the limits $n={\rm const}\gg 1$, $\mu\gg 1$: 
\be \label{V1}
V^{(1)}(n,\mu)=\frac{32\pi^2}{\ln^42}\;\frac{2^n}{n^3}\;7^{\mu-1} 
\ee
\bigskip

\noindent {\sc The group ${\cal LF}_{n+1}^{(2)}$}. It is convenient to 
enumerate the generators $f_{ij}^{(\alpha)}$, $i,j=1,\ldots,n^2$, 
$\alpha=x,y$, ordering them in a sequence: $(f_{11}^{(x)}, f_{11}^{(y)}), 
\ldots, (f_{1n}^{(x)}, f_{1n}^{(y)}), \ldots, (f_{n1}^{(x)}, f_{n1}^{(y)}), 
\ldots, (f_{nn}^{(x)}, f_{nn}^{(y)})$. For any such sequence we define the 
"normal order" according to the prescriptions (i)-(iii). Let $z$ be the 
serial number of the pair $(f_{ij}^{(x)},f_{ij}^{(y)})$. Consider the 
functions $a(z,m)$ and $b(z,m)$ defined as numbers of all distinct 
sequences of $m+1$ generators satisfying the rules (i)-(iii) and ending 
with $f_{ij}^{(x)}$ and $f_{ij}^{(y)}$ respectively. One can readilly write 
the evolution equations for $a(z,m)$ and $b(z,m)$ similar to (\ref{5}):
\be \left\{\begin{array}{lll}
a(z,m+1) & = & a(z-1,m)+b(z-n,m)+b(z-n+1,m)+b(z,m) +  \\
& & \disp \sum_{z'=z+1}^{n^2} \Big( a(z',m)+b(z',m) \Big) \medskip \\ 
b(z,m+1) & = & \disp b(z-n,m)+a(z-1,m)+a(z,m)+ \sum_{z'=z+1}^{n^2} 
\Big( a(z',m)+b(z',m) \Big) \end{array} \right.  
\ee

Analogous to the case of the group ${\cal LF}_{n+1}^{(1)}$ let us suppose 
in the limit $n \gg 1$ the periodical boundary conditions on the segment 
$[0,n^2+1]$. So we have
\be \label{ab} 
\left\{\begin{array}{l}
\begin{array}{lll}
a(z+1,m+1)-a(z,m+1) & = & b(z-n+2,m)-b(z-n,m)-b(z,m)+ \\
& & a(z,m)-a(z-1,m)-a(z+1,m) \\
b(z+1,m+1)-b(z,m+1) & = & b(z-n+1,m)-b(z-n,m)-a(z-1,m)- \\
& & b(z+1,m) 
\end{array} \medskip \\
~a(0,m)=b(0,m)=a(n^2+1,m)=b(n^2+1,m)=0, \qquad m=0,1,\ldots \medskip \\
~a(z,0)=b(z,0)=1, \qquad z=1,\ldots,n^2.
\end{array} \right.
\ee 
The initial conditions differ from (\ref{theta}) because in (\ref{ab}) we 
do not fix the first generator in the sequence involved. 
 
Assuming that $|a(z,m)-b(z,m)|\to 0$ ($n\to\infty$) uniformly for $z$ and 
$m$, we may pass from (\ref{ab}) to a single closed equation for the 
function $b(z,m)$. (The selfconsistency of this supposition we check at the 
end of our computations). So, we get:
\be 
\left\{\begin{array}{l}
\begin{array}{lll}
b(z+1,m+1)-b(z,m+1) & = & b(z-n+1,m)-b(z-n,m)-b(z-1,m)- \\ 
& & b(z+1,m) 
\end{array} \medskip \\
~b(0,m)=b(n^2+1,m)=0, \qquad m=0,1, \ldots \medskip \\  
~b(z,0)=1, \qquad z=1,\ldots,n^2.
\end{array} \right.
\ee
Performing the decomposition $\disp b(z,m)=\sum_{k=1}^{n^2} 
B_k\lam_k^m\beta_k(z)$, we arrive at the following boundary problem:  
\be \label{beta} 
\left\{\begin{array}{l}
(\lam_k+1)\beta_k(z+1)-\lam_k\beta_k(z)+\beta_k(z-1)+\beta_k(z-n)-
\beta_k(z-n+1)=0 \\
\beta_k(0)=\beta_k(n^2+1)=0.
\end{array} \right.
\ee

Let us look for the solution of Eq.(\ref{beta}) in the form $\disp \beta_k(z)= 
p_k^z\sin\frac{\pi k z}{n^2+1}$. Substituting this ansatz in (\ref{beta}) 
we obtain an equation for $p_k$ as well as an expression for $\lam_k$:
\be \label{p}
\left\{\begin{array}{l}
\begin{array}{r}
\disp \sin\frac{\pi k}{n^2+1}p_k^{n+1}+ \sin\frac{2\pi k}{n^2+1}p_k^n-
\sin\frac{\pi k}{n^2+1}p_k^{n-1}- \sin\frac{\pi kn}{n^2+1}p_k^2+ \\
\disp 2\cos\frac{\pi k}{n^2+1}\sin\frac{\pi kn}{n^2+1}p_k- \sin\frac{\pi 
kn}{n^2+1}=0 
\end{array} \medskip \\ 
\disp ~\lam_k = p_k^{-2}-1-p_k^{-n}\frac{\sin\frac{\pi 
k(n-1)}{n^2+1}}{\sin\frac{\pi k}{n^2+1}}+p_k^{-n-1}\frac{\sin\frac{\pi 
kn}{n^2+1}}{\sin\frac{\pi k}{n^2+1}} 
\end{array}\right. 
\ee

In Eq.(\ref{p}) each root $p_k^{(i)}$ corresponds to different values of
$\lam_k^{(i)}$. However, we are interested only in the asymptotic behavior of 
$b(z,m)$ ($m\gg 1$) determined by the largest value of $\lam_k^{(i)}$ 
for $k=1$. (Compare to the case of the group ${\cal 
LF}_{n+1}^{(1)}$---Eq.(\ref{eigen1d})). The Eq.(\ref{p}) at $k=1$ and $n 
\gg 1$ reads:
\be \label{p1}
\left\{\begin{array}{l}
p_1^{n+1}+ 2p_1^n-p_1^{n-1}-n(p_1-1)^2=0 \medskip \\
\lam_1=p_1^{-2}-1+np_1^{-n-1}-(n-1)p_1^{-n}.
\end{array}\right.
\ee
One can easily check that the smallest positive root corresponding to the 
largest value of $\lam_1^{(i)}$ is 
$$
p_1=1-\frac{\ln n}{n}+o\left( \frac{\ln n}{n} \right), \qquad 
n={\rm const}\gg 1 
$$ 
and
\be \label{lambda1} 
\lam_1=\frac{n}{\ln n} + o\left(\frac{n}{\ln n} \right), \qquad 
n={\rm const}\gg 1
\ee 
In the $1+1$--dimensional case the same value of $\lam_k$ was given 
by the right edge of the spectrum, but in $2+1$--dimensional case one can
prove that Eq.(\ref{p}) for any $i>1$ has no solutions $\lam_{n^2}^{(i)}$ 
growing as fast as $\lam_1$. The coefficients $B_k$ should be found from 
the initial condition 
$$
\sum_{k=1}^{n^2}B_k\beta_k(z)=1.
$$ 
As $p_1 \to 1$ ($n\gg 1$) we can except the set $\beta_k(z)$ to be 
orthogonal in the vicinity of the left edge of the spectrum, so $B_1$ is 
determined basically by the expression 
$$
B_1=\frac{\disp\int_0^{n^2+1}\beta_1(z)dz}
{\disp \int_0^{n^2+1}\beta_1^2(z)dz}, \quad n\gg 1 
$$ 
Now we have the following equation for the function for $b(z,m)$ in the 
limits $n={\rm const} \gg 1$ and $m \gg 1$ 
\be \label{b}
b(z,m)=\frac{4}{\pi}\sin\frac{\pi z}{n^2+1} 
\left(\frac{n}{\ln n}\right)^m
\ee

If we suppose the equality $a(z,m)=b(z,m)$, where $b(z,m)$ is given by 
(\ref{b}), it is easy to check that $a(z,m)$ and $b(z,m)$ really satisfy 
the equations (\ref{ab}) in the limits $n \gg 1$, $m \gg 1$. This 
fact proves our assumption about the behaviors of $a(z,m)$ and $b(z,m)$ in 
a selfconsistent way.

The limiting expression of the function $\tilde{\th}_n(m)$ reads
$$
\tilde{\th}_n(m)=\sum\limits_{z=1}^{n^2} a(z,m)+b(z,m)=
\frac{16n^2}{\pi ^2}\left(
\frac{n}{\ln n}\right)^m.
$$ 
Thus, the asymptotics of the number of nonequivalent words of given 
irreducible length, $\mu$ in the limit $n={\rm const} \gg 1$, $\mu\gg 1$ is
\be \label{V2} 
V^{(2)}(n,\mu)=\sum_{m=1}^{\mu-1}2^{m+1}\left(\mu-1 \atop m\right)
\tilde{\th}_n(m)=\frac{32n^2}{\pi^2}\left(\frac{2n}{\ln n}\right)^{\mu-1} 
\ee  
(compare to Eq.(\ref{V1})).

\section{Conclusion}

The principal difference between the limitng behavior of the partition 
functions $V^{(1)}(n,\mu)$ and $V^{(2)}(n,\mu)$ (and, hence, between the
upper boundaries of the sets $\Omega^{(1)}(n,\mu)$ and 
$\Omega^{(2)}(n,\mu)$) becomes at most illuminating in the limit 
$n={\rm const}\gg 1$ and $\mu\to \infty$ if we consider the following 
limit 
\be \label{f}
f_{1,2}=\left[\lim_{\mu\to\infty}\frac{\ln V^{(1,2)}(n,\mu)}{\mu}
\right]_{n={\rm const}\gg 1}
\ee

Using the equations (\ref{V1}) and (\ref{V2}), we get
$$
\left\{\begin{array}{l}
f_1=\ln 7 \medskip \\ \disp f_2=\frac{2n}{\ln n}
\end{array}\right.
$$

Thus, we can conclude, that with the exponential accuracy in the limit 
$n={\rm const}\gg 1$ and $\mu\to \infty$ the set $\Omega^{(1)}(n,\mu)$
is bounded from above by the $n$--independent estimate, i.e. 
$\Omega^{(1)}(n,\mu)$ is "representation--independent"; while the set 
$\Omega^{(2)}(n,\mu)$ with the same accuracy and in the same limit depends
strongly on the braid representation (i.e. on the number of strings, $n$).

The equations (\ref{V1}) and (\ref{V2}) enable us to make some conclusions 
about the structure of the graphs corresponding to the groups 
${\cal LF}_n^{(1)}$ and ${\cal LF}_n^{(2)}$. These graphs can be viewed as 
follows. Take the {\it free} $1+1$-- or $2+1$-- groups, where all generators 
do not commute at all. The graphs of these groups have structures of $2n$-- 
and $4n^2$-- branching Cayley trees, where the number of distinct words of 
length $\mu$ is equal to 
$$
\left\{\begin{array}{ll}
V_{\rm free}^{(1)}(n,\mu)=2n(2n-1)^{\mu-1} & \mbox{for $1+1$-- free group} 
\medskip \\
V_{\rm free}^{(2)}(n,\mu)=4n^2(4n^2-1)^{\mu-1} & \mbox{for $2+1$-- free group}
\end{array}\right.
$$

The graphs corresponding to the groups ${\cal LF}_n^{(1,2)}$
can be constructed from the graphs of the free groups in accordance with the
following recursion procedure: 
\begin{itemize}
\item[(i)] Take the root vertex of the free group graph and consider all 
vertices on the distance $\mu=2$. Identify those vertices which correspond to 
the equivalent words in groups ${\cal LF}_n^{(1,2)}$; 
\item[(ii)] Repeat this procedure taking all vertices on the distance 
$\mu=(1,2,\ldots)$ and ``gluing" them on the distance $\mu+2$ according 
to the definition of the locally free groups. 
\end{itemize}
By means of this procedure we raise a graph which in average has 
$z_{\rm eff}^{(1,2)}-1$ distinct branches leading from the level $\mu$ to the 
level $\mu+1$. We may easily find the expressions of $z_{\rm eff}^{(1,2)}$ 
using the Eqs. (\ref{V1}) and (\ref{V2}). We have in the limit 
$n={\rm const}\gg 1$ and $\mu\to \infty$:
$$
\disp z_{\rm eff}^{(1,2)}=\frac{V^{(1,2)}(n,\mu+1)}{V^{(1,2)}(n,\mu)}+1 =
\left\{\begin{array}{ll}
8 & \quad \mbox{for $1+1$-- locally free group} \medskip \\
\disp \frac{2n}{\ln n}+1 & \quad \mbox{for $2+1$-- locally free group}
\end{array}\right.
$$
We see that the graph of the group ${\cal LF}_n^{(1)}$ coincides (in average) 
with $(z_{\rm eff}=8)$--branching Cayley tree for any $n\gg 1$, while the 
effective coordinational number of the graph of the group ${\cal LF}_n^{(2)}$ 
depends on $n$ and does not "saturare" for $n\gg 1$.

\noindent {\bf Acknowledgments}

We are very grateful to A. Comtet, J. Desbois and A.M. Vershik for fruitful
discussions and usful suggestions.

\newpage
\section*{Figure Captions}

\begin{fig}
Representation of: (a) $1+1$--braid group generator, $\sigma_i$; (b) 
$2+1$--braid group generators $\sigma_{i,j}^{x}$ and $\sigma_{i,j}^{y}$ 
\label{fig1} 
\end{fig}

\begin{fig}
Representation of: (a) $1+1$--locally free group generator, $f_i$; (b) 
$2+1$--locally free group generators $f_{i,j}^{x}$ and $f_{i,j}^{y}$ 
\label{fig2} 
\end{fig}


\newpage

\begin{thebibliography}{99}

\bibitem{pr} L. Paris, D. Rolfsen, {\it Geometric subgroups of surface 
braid groups}, Preprint No. 115 (1997) (Universit\'e de Bourgogne, 
Laboratoire de Topologie)

\bibitem{birm} J. Birman, {\it Mapping class groups of surfaces}, Contemp. 
Math., 78 (1988), 13 

\bibitem{ver1} A.M. Vershik, S.V. Kerov, Sov. Ac. Sci. Doklady, 301 (1988), 
777; A.M. Vershik, Topics in Algebra, 26 (1990), pt. 2, 467; 
Proc. Am. Math. Soc. 148 (1991), 1

\bibitem{negrve} S.K. Nechaev, A.Yu. Grosberg, A.M. Vershik, J. Phys. A: 
Math. Gen., 29 (1996), 2411; J. Desbois, S. Nechaev, J. Stat. Phys., 88 
(1997), 201-223; J. Desbois, S. Nechaev, J. Phys. (A): Math. Gen., 31 (1998), 
2767-2784;  A. Comtet, S. Nechaev, J. Phys. (A): Math. Gen. 31 (1998),

\bibitem{cafo} P. Cartier, D. Foata, Lect. Notes in Math., 85 (1969)

\bibitem{vien} G.X. Viennot, Lect. Notes in Math., 1234 (1986), 321

\bibitem{ver2} A.M. Vershik, privite communication 

\end{thebibliography}
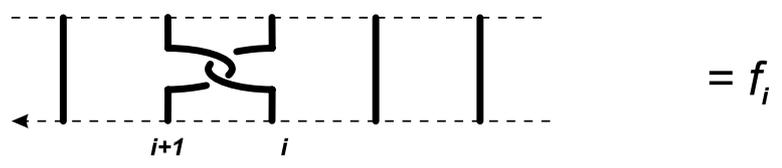
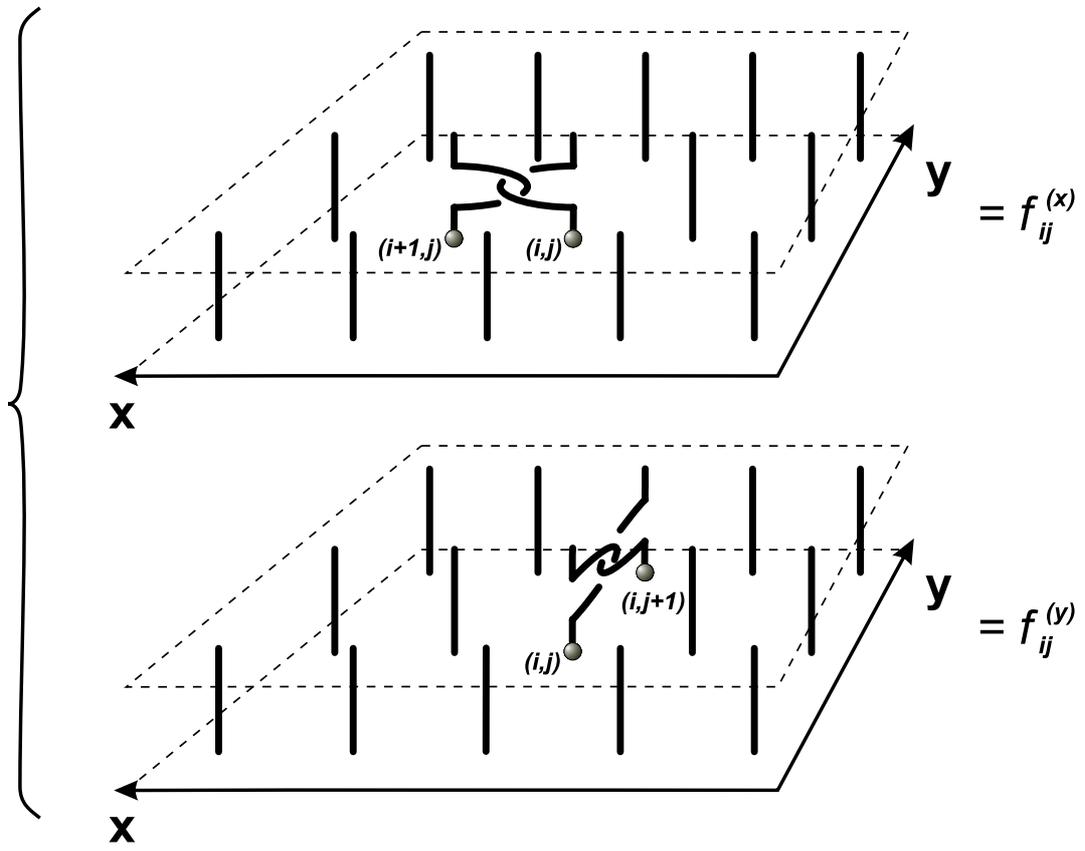
\end{document}